\documentclass[a4paper,11pt]{article}
\usepackage{amsmath,amsthm,amssymb}
\usepackage[mathscr]{eucal}

{\theoremstyle{definition}\newtheorem{definition}{Definition}[section]

\newtheorem{remark}[definition]{Remark}
}

\newtheorem{proposition}[definition]{Proposition}
\newtheorem{lemma}[definition]{Lemma}
\newtheorem{theorem}[definition]{Theorem}
\newtheorem{corollary}[definition]{Corollary}

\newcommand{\rc}{\operatorname{c}}
\newcommand{\rC}{\operatorname{C}}
\newcommand{\Cred}{C_\text{\rm red}}
\newcommand{\otmin}{\ot_{\text{\rm min}}}
\newcommand{\otalg}{\ot_{\text{\rm alg}}}
\newcommand{\mult}{\operatorname{mult}}
\newcommand{\Tr}{\operatorname{Tr}}
\newcommand{\dK}{\mathbb{K}}
\newcommand{\U}{\operatorname{U}}
\newcommand{\SU}{\operatorname{SU}}

\newcommand{\C}{\mathbb{C}}
\newcommand{\GL}{\operatorname{GL}}
\newcommand{\Hcentr}{\operatorname{H^\infty_\text{\rm centr}}}
\newcommand{\Irred}{\operatorname{Irred}}
\newcommand{\Ker}{\operatorname{Ker}}
\newcommand{\Mor}{\operatorname{Mor}}
\newcommand{\M}{\operatorname{M}}
\newcommand{\N}{\mathbb{N}}

\newcommand{\PI}{\Theta_{\mu}}

\newcommand{\R}{\mathbb{R}}
\newcommand{\al}{\alpha}
\newcommand{\atil}{\widetilde{a}}
\newcommand{\bP}{\mathbb{P}}
\newcommand{\bV}{\mathbb{V}}
\newcommand{\be}{\beta}
\newcommand{\bnd}{\operatorname{bnd}}

\newcommand{\cE}{\mathcal{E}}

\newcommand{\cGh}{\widehat{\mathbb{G}}}
\newcommand{\cG}{\mathbb{G}}

\newcommand{\cK}{\mathcal{K}}
\newcommand{\cL}{\mathcal{L}}

\newcommand{\cU}{\mathcal{U}}
\newcommand{\cV}{\mathcal{V}}

\newcommand{\cZ}{\mathcal{Z}}
\newcommand{\deh}{\hat{\Delta}}
\newcommand{\de}{\Delta}
\newcommand{\dimq}{\operatorname{dim}_q}
\newcommand{\dpr}{^{\prime\prime}}
\newcommand{\epsh}{\widehat{\epsilon}}
\newcommand{\eps}{\epsilon}
\newcommand{\fancyot}{\mathbin{\text{\footnotesize\textcircled{\tiny \sf T}}}}
\newcommand{\id}{\mathord{\text{\rm id}}}

\newcommand{\lx}{\vert x\vert}

\newcommand{\mB}{\mathscr{B}}
\newcommand{\om}{\omega}
\newcommand{\ot}{\otimes}
\newcommand{\ovt}{\overline{\otimes}}

\newcommand{\poisson}{\operatorname{H^{\infty}}(\widehat{\mathbb{G}},\mu)}
\newcommand{\rL}{\operatorname{L}}

\newcommand{\recht}{\rightarrow}
\newcommand{\si}{\sigma}
\newcommand{\supp}{\operatorname{supp}}
\newcommand{\ux}{\underline{x}}
\newcommand{\veps}{\varepsilon}

\newcommand{\lspan}{\operatorname{span}}
\newcommand{\dimmin}{\operatorname{dim}_{\text{\rm min}}}

\begin{document}
\begin{center}
{\LARGE\bf Poisson boundary of the discrete quantum group $\widehat{A_u(F)}$}

\bigskip

{\sc by Stefaan Vaes and Nikolas Vander Vennet\renewcommand{\thefootnote}{}\footnotetext{\parbox[t]{450pt}{Department of Mathematics;
K.U.Leuven; Celestijnenlaan 200B; B--3001 Leuven (Belgium)\\ E-mail: stefaan.vaes@wis.kuleuven.be and nikolas.vandervennet@wis.kuleuven.be \\ Research partially supported by ERC Starting Grant VNALG-200749 and Research
    Programme G.0231.07 of the Research Foundation --
    Flanders (FWO)}}}
\end{center}

\begin{abstract}
\noindent We identify the Poisson boundary of the dual of the universal compact quantum group $A_u(F)$ with a measurable field of ITPFI factors.
\end{abstract}

\section{Introduction and statement of main result}

Poisson boundaries of discrete quantum groups were introduced by Izumi \cite{Iz} in his study of infinite tensor product actions of $\SU_q(2)$. Izumi was able to identify the Poisson boundary of the dual of $\SU_q(2)$ with the quantum homogeneous space $\rL^\infty(\SU_q(2)/S^1)$, called the Podle{\'s} sphere. The generalization to $\SU_q(n)$ was established by Izumi, Neshveyev and Tuset \cite{INT}, yielding $\rL^\infty(\SU_q(n)/S^{n-1})$ as the Poisson boundary. A more systematic approach was given by Tomatsu \cite{Tom} who proved the following very general result: if $\cG$ is a compact quantum group with commutative fusion rules and amenable dual $\cGh$, the Poisson boundary of $\cGh$ can be identified with the quantum homogeneous space $\rL^\infty(\cG/\dK)$, where $\dK$ is the maximal closed quantum subgroup of Kac type inside $\cG$. Tomatsu's result provides the Poisson boundary for the duals of all $q$-deformations of classical compact groups.

All examples discussed in the previous paragraph concern amenable discrete quantum groups. In \cite{VaVV}, we identified the Poisson boundary for the (non-amenable) dual of the compact quantum group $A_o(F)$ with a higher dimensional Podle{\'s} sphere. Although the dual of $A_o(F)$ is non-amenable, the representation category of $A_o(F)$ is monoidally equivalent with the representation category of $\SU_q(2)$ for the appropriate value of $q$. The second author and De Rijdt provided in \cite{DRVV} a general result explaining the behavior of the Poisson boundary under the passage to monoidally equivalent quantum groups. In particular, a combination of the results of \cite{DRVV} and \cite{Iz} give a more conceptual approach to our identification in \cite{VaVV}.

The quantum random walks studied on a discrete quantum group $\cGh$ have a semi-classical counterpart, being a Markov chain on the (countable) set $\Irred(\cG)$ of irreducible representations of $\cG$ (modulo unitary equivalence). All the examples above, share the feature that the semi-classical random walk on $\Irred(\cG)$ has trivial Poisson boundary.

In this paper, we identify the Poisson boundary for the dual of $\cG = A_u(F)$. In that case, $\Irred(\cG)$ can be identified with the Cayley tree of the monoid $\N * \N$ and, by results of \cite{PiWo}, has a non-trivial Poisson boundary: the end compactification of the tree with the appropriate harmonic measure. Before discussing in more detail our main result, we introduce some terminology and notations. For a more complete introduction to Poisson boundaries of discrete quantum groups, we refer to \cite[Chapter 4]{VV-thesis}.

Compact quantum groups were originally introduced by Woronowicz in \cite{Wor1} and their definition finally took the following form.

\begin{definition}[{Woronowicz \cite[Definition 1.1]{Wor2}}]
A \emph{compact quantum group} $\cG$ is a pair consisting of a unital C$^*$-algebra $\rC(\cG)$ and a unital $*$-homomorphism $\Delta : \rC(\cG) \recht \rC(\cG) \ot \rC(\cG)$, called \emph{comultiplication,} satisfying the following two conditions.
\begin{itemize}
\item \emph{Co-associativity:} $(\de \ot \id)\de = (\id \ot \de) \de$.
\item $\lspan \de(\rC(\cG)) (1 \ot \rC(\cG))$ and $\lspan \de(\rC(\cG)) (\rC(\cG) \ot 1)$ are dense in $\rC(\cG) \ot \rC(\cG)$.
\end{itemize}
\end{definition}
In the above definition, the symbol $\ot$ denotes the minimal (i.e.\ spatial) tensor product of C$^*$-algebras.

Let $\cG$ be a compact quantum group. By \cite[Theorem 1.3]{Wor2}, there is a unique state $h$ on $\rC(\cG)$ satisfying $(\id \ot h)\de(a) = h(a)1 = (h \ot \id)\de(a)$ for all $a \in \rC(\cG)$. We call $h$ the \emph{Haar state} of $\cG$.

A \emph{unitary representation} of $\cG$ on the finite dimensional Hilbert space $H$ is a unitary operator $U \in \cL(H) \ot \rC(\cG)$ satisfying $(\id \ot \de)(U) = U_{12} U_{13}$. Given unitary representations $U_1,U_2$ on $H_1,H_2$, we put
$$\Mor(U_2,U_1) := \{S \in \cL(H_1,H_2) \mid (S \ot 1) U_1 = U_2(S \ot 1) \} \; .$$
Let $U$ be a unitary representation of $\cG$ on the finite dimensional Hilbert space $H$. The elements $(\xi^* \ot 1) U (\eta \ot 1) \in \rC(\cG)$ are called the coefficients of $U$. The linear span of all coefficients of all finite dimensional unitary representations of $\cG$ forms a dense $*$-subalgebra of $\rC(\cG)$ (see \cite[Theorem 1.2]{Wor2}). We call $U$ \emph{irreducible} if $\Mor(U,U) = \C 1$. We call $U_1$ and $U_2$ \emph{unitarily equivalent} if $\Mor(U_2,U_1)$ contains a unitary operator.

Let $U$ be an irreducible unitary representation of $\cG$ on the finite dimensional Hilbert space $H$. By \cite[Proposition 5.2]{Wor2}, there exists an anti-linear invertible map $j : H \recht \overline{H}$ such that the operator $U^c \in \cL(\overline{H}) \ot \rC(\cG)$ defined by the formula $(j(\xi)^* \ot 1)U^c(j(\eta) \ot 1) = (\eta^* \ot 1)U^* (\xi \ot 1)$ is unitary. One calls $U^c$ the \emph{contragredient} of $U$. Since $U$ is irreducible, the map $j$ is uniquely determined up to multiplication by a non-zero scalar. We normalize in such a way that $Q := j^* j$ satisfies $\Tr(Q) = \Tr(Q^{-1})$. Then, $j$ is determined up to multiplication by $\lambda \in S^1$ and $Q$ is uniquely determined. We call $\Tr(Q)$ the \emph{quantum dimension} of $U$ and denote it by $\dimq(U)$. Note that $\dimq(U) \geq \dim(H)$ with equality holding iff $Q = 1$.

The \emph{tensor product} $U \fancyot V$ of two unitary representations is defined as $U_{13} V_{23}$.

Given a compact quantum group $\cG$, we denote by $\Irred(\cG)$ the set of irreducible unitary representations of $\cG$ modulo unitary conjugacy. For every $x \in \Irred(\cG)$, we choose a representative $U^x$ on the Hilbert space $H_x$. We denote by $Q_x \in \cL(H_x)$ the associated positive invertible operator and define the state $\psi_x$ on $\cL(H_x)$ by the formula
$$\psi_x(A) := \frac{\Tr(Q_x A)}{\Tr(Q_x)} \; .$$
The dual, \emph{discrete quantum group} $\cGh$ is defined as the $\ell^\infty$-direct sum of matrix algebras
$$\ell^\infty(\cGh) := \prod_{x \in \Irred(\cG)} \cL(H_x) \; .$$
We denote by $p_x, x \in \Irred(\cG)$, the minimal central projections in $\ell^\infty(\cGh)$. Denote by $\eps \in \Irred(\cG)$ the trivial representation and by $\epsh : \ell^\infty(\cGh) \recht \C$ the co-unit given by $a p_\eps = \epsh(a) p_\eps$.

Whenever $x,y,z \in I$, we use the short-hand notation $\Mor(x \ot y,z) := \Mor(U^x \fancyot U^y , U^z)$ and we write $z \subset x \ot y$ if $\Mor(x \ot y,z) \neq \{0\}$.

The von Neumann algebra $\ell^\infty(\cGh)$ carries a comultiplication $\deh : \ell^\infty(\cGh) \recht \ell^\infty(\cGh) \ovt \ell^\infty(\cGh)$, uniquely characterized by the formula
$$\deh(a) (p_x \ot p_y) S = S a p_z \quad\text{for all}\quad x,y,z \in \Irred(\cG) \;\;\text{and}\;\; S \in \Mor(x \ot y,z) \; .$$
Denote by $\rL^\infty(\cG)$ the weak closure of $\rC(\cG)$ in the GNS representation of the Haar state $h$. One defines the unitary $\bV \in \ell^\infty(\cGh) \ovt \rL^\infty(\cG)$ by the formula
$$\bV := \bigoplus_{x \in \Irred(\cG)} U^x \; .$$
The unitary $\bV$ implements the duality between $\cG$ and $\cGh$, in the sense that it satisfies
$$(\deh \ot \id)(\bV) = \bV_{13} \bV_{23} \quad\text{and}\quad (\id \ot \de)(\bV) = \bV_{12} \bV_{13} \; .$$
Discrete quantum groups can also be defined intrinsically, see \cite{VD1}.

Whenever $\om \in \ell^\infty(\cGh)_*$ is a normal state, we consider the Markov operator
$$P_\om : \ell^\infty(\cGh) \recht \ell^\infty(\cGh) : P_\om(a) = (\id \ot \om)\deh(a) \; .$$
By \cite[Proposition 2.1]{NT}, the Markov operator $P_\om$ leaves globally invariant the center $\cZ(\ell^\infty(\cGh))$ of $\ell^\infty(\cGh)$ if and only if
$$\om = \psi_\mu := \sum_{x \in \Irred(\cG)} \mu(x) \psi_x \quad\text{where}\quad \mu \;\;\text{is a probability measure on}\;\; \Irred(\cG) \; .$$
We only consider states $\om$ of the form $\psi_\mu$ and denote by $P_\mu$ the corresponding Markov operator. Note that we can define a convolution product on the probability measures on $\Irred(\cG)$ by the formula
$$P_{\mu * \eta} = P_\mu \circ P_\eta \; .$$
Considering the restriction of $P_\mu$ to $\ell^\infty(\Irred(\cGh)) = \cZ(\ell^\infty(\cGh))$,
every probability measure $\mu$ on $\Irred(\cG)$ defines a Markov chain on the countable set $\Irred(\cG)$ with $n$-step transition probabilities given by
$$p_x p_n(x,y) = p_x P_\mu^n(p_y) \; .$$
The probability measure $\mu$ is called \emph{generating} if for every $x,y \in \Irred(\cG)$, there exists an $n \in \N \setminus \{0\}$ such that $p_n(x,y) > 0$.

\begin{definition} \label{def.poisson}
Let $\cG$ be a compact quantum group and $\mu$ a generating probability measure on $\Irred(\cG)$. The \emph{Poisson boundary} of $\cGh$ with respect to $\mu$ is defined as the space of $P_\mu$-harmonic elements in $\ell^\infty(\cGh)$.
$$\poisson := \{a \in \ell^\infty(\cGh) \mid P_\mu(a) = a \} \; .$$
The weakly closed vector subspace $\poisson$ of $\ell^\infty(\cGh)$ is turned into a von Neumann algebra using the product (cf.\ \cite[Theorem 3.6]{Iz})
$$a \cdot b := \lim_{n \recht \infty} P_\mu^n(ab)$$
and where the sequence at the right hand side is strongly$^*$ convergent.
\begin{itemize}
\item The restriction of $\epsh$ to $\poisson$ is a faithful normal state on $\poisson$.
\item The restriction of $\deh$ to $\poisson$ defines a left action
$$\al_{\cGh} : \poisson \recht \ell^\infty(\cGh) \ovt \poisson : a \mapsto \deh(a)$$
of $\cGh$ on $\poisson$.
\item The restriction of the adjoint action to $\poisson$ defines an action
$$\al_\cG : \poisson \recht \poisson \ovt \rL^\infty(\cG) : a \mapsto \bV (a \ot 1) \bV \; .$$
\end{itemize}
\end{definition}

We denote by $\Hcentr(\cGh,\mu) := \poisson \cap \cZ(\ell^\infty(\cGh))$ the space of bounded $P_\mu$-harmonic functions on $\Irred(\cG)$. Defining the conditional expectation
$$\cE : \ell^\infty(\cGh) \recht \ell^\infty(\Irred(\cGh)) : \cE(a) p_x = \psi_x(a) p_x \; ,$$
we observe that $\cE$ also provides a faithful conditional expectation of $\poisson$ onto the von Neumann subalgebra $\Hcentr(\cGh,\mu)$.

We now turn to the concrete family of compact quantum groups studied in this paper and introduced by Van Daele and Wang in \cite{VDW}. Let $n \in \N \setminus \{0,1\}$ and let $F \in \GL(n,\C)$. One defines the compact quantum group $\cG = A_u(F)$ such that $\rC(\cG)$ is the universal unital C$^*$-algebra generated by the entries of an $n \times n$ matrix $U$ satisfying the relations
$$U \quad\text{and}\quad F \overline{U} F^{-1} \quad\text{are unitary, with}\quad \bigl(\overline{U}\bigr)_{ij} = (U_{ij})^*$$
and such that $\de(U_{ij}) = \sum_{k=1}^n U_{ik} \ot U_{kj}$. By definition, $U$ is an $n$-dimensional unitary representation of $A_u(F)$, called the \emph{fundamental representation.}

Fix $F \in \GL(n,\C)$ and put $\cG = A_u(F)$. For reasons to become clear later, we assume that $F$ is not a scalar multiple of a unitary $2 \times 2$ matrix.

By \cite[Th\'{e}or\`{e}me 1]{Ba2}, the irreducible unitary representations of $\cG$ can be labeled by the elements of the free monoid $I := \N * \N$ generated by $\alpha$ and $\beta$. We represent the elements of $I$ as words in $\alpha$ and $\beta$. The empty word is denoted by $\eps$ and corresponds to the trivial representation of $\cG$, while $\al$ corresponds to the fundamental representation and $\be$ to the contragredient of $\al$. We denote by $x \mapsto \overline{x}$ the unique antimultiplicative and involutive map on $I$ satisfying $\overline{\alpha} = \beta$. This involution corresponds to the contragredient on the level of representations. The fusion rules of $\cG$ are given by
$$x \ot y \cong \bigoplus_{z \in I, x = x_0 z, y = \overline{z} y_0} x_0 y_0 \; .$$
So, if the last letter of $x$ equals the first letter of $y$, the tensor product $x \ot y$ is irreducible and given by $xy$. We denote this as $xy = x \ot y$.

Denote by $\partial I$ the compact space of infinite words in $\alpha$ and $\beta$. For $x \in \partial I$, the expression
\begin{equation}\label{eq.infiniteproduct}
x = x_1 \otimes x_2 \otimes \cdots
\end{equation}
means that the infinite word $x$ is the concatenation of the finite words $x_1x_2\cdots$ and that the last letter of $x_n$ equals the first letter of $x_{n+1}$ for all $n \in \N$. All elements $x$ of $\partial I$ can be decomposed as in \eqref{eq.infiniteproduct}, except the countable number of elements of the form
$x = y \alpha \beta \alpha \beta \cdots$ for some $y \in I$.

Below, we will only deal with non-atomic measures on $\partial I$, so that almost every point of $\partial I$ has a decomposition as in \eqref{eq.infiniteproduct}. We denote by $\partial_0 I$ the subset of $\partial I$ consisting of the infinite words that have a decomposition of the form \eqref{eq.infiniteproduct}.

The following is the main result of the paper.

\begin{theorem}\label{thm.main}
Let $F \in \GL(n,\C)$ such that $F$ is not a scalar multiple of a unitary $2 \times 2$ matrix. Write $\cG = A_u(F)$ and suppose that $\mu$ is a finitely supported, generating probability measure on $I=\Irred(\cG)$. Denote by $\partial I$ the compact space of infinite words in the letters $\alpha,\beta$. There exists
\begin{itemize}
\item a non-atomic probability measure $\nu_\eps$ on $\partial I$,
\item a measurable field $M$ of ITPFI factors over $(\partial I,\nu_\eps)$ with fibers
$$(M_x,\omega_x) = \bigotimes_{k=1}^\infty (\cL(H_{x_k}),\psi_{x_k})$$
whenever $x \in \partial_0 I$ is of the form $x = x_1x_2x_3 \cdots = x_1 \otimes x_2 \otimes x_3 \otimes \cdots$,
\item an action $\beta_{\cGh}$ of $\cGh$ on $M$ concretely given by \eqref{eq.concrete} below,
\end{itemize}
such that, with $\om_\infty = \int^{\oplus} \omega_x \; d\nu_\eps(x)$, the Poisson integral formula
$$\PI : M \recht \poisson : \PI(a) = (\id \otimes \om_\infty)\beta_{\cGh}(a)$$
defines a $*$-isomorphism of $M$ onto $\poisson$, intertwining the action $\beta_{\cGh}$ on $M$ with the action $\alpha_{\cGh}$ on $\poisson$.

Moreover, defining the action $\beta^x_{\cG}$ of $\cG$ on $M_x$ as the infinite tensor product of the inner actions $a \mapsto  U^{x_k}(a\ot 1)(U^{x_k})^*$, we obtain the action $\beta_\cG$ of $\cG$ on $M$. The $*$-isomorphism $\PI$ intertwines $\beta_\cG$ with $\alpha_\cG$.
\end{theorem}

The comultiplication $\deh : \ell^\infty(\cGh) \recht \ell^\infty(\cGh) \ovt \ell^\infty(\cGh)$ can be uniquely cut down into completely positive maps $\deh_{x \ot y,z} : \cL(H_z) \recht \cL(H_x) \ot \cL(H_y)$ in such a way that
$$\deh(a)(p_x \ot p_y) = \sum_{z \subset x \ot y} \deh_{x \ot y,z}(ap_z)$$
for all $a \in \ell^\infty(\cGh)$.

We denote by $|x|$ the length of a word $x \in I$.

If now $x,y \in I, z \in \partial I$ with $yz = y \otimes z$ and $|y|>|x|$, we define for all $s \subset x \otimes y$,
$$\deh_{x \ot yz,sz} : M_{sz} \recht \cL(H_x) \ot M_{yz}$$
by composing $\deh_{x \ot y,s} \ot \id$ with the identifications $M_{sz} \cong \cL(H_s) \ot M_z$ and $M_{yz} \cong \cL(H_y) \ot M_z$. The action $\beta_{\cGh} : M \recht \ell^\infty(\cGh) \ovt M$ of $\cGh$ on $M$ is now given by
\begin{equation}\label{eq.concrete}
\beta_{\cGh}(a)_{x,yz} = \sum_{s \subset x \ot y} \deh_{x \ot yz,sz}(a p_{sz})
\end{equation}
whenever $a \in M$, $x,y \in I$, $z \in \partial I$, $|y| > |x|$ and $yz = y \ot z$. Note that we identified $\ell^\infty(\cGh) \ovt M$ with a measurable field over $I \times \partial I$ with fiber in $(x,z)$ given by $\cL(H_x) \ot M_z$.

\subsection*{Further notations and terminology}

Fix $F \in \GL(n,\C)$ and put $\cG = A_u(F)$. We identify $\Irred(\cG)$ with $I := \N * \N$. We assume that $F$ is not a multiple of a unitary $2 \times 2$ matrix. Equivalently, $\dimq(\al) > 2$. The first reason to do so, is that under this assumption, the random walk defined by any non-trivial probability measure $\mu$ on $I$ (i.e.\ $\mu(\eps) < 1$), is automatically \emph{transient,} which means that
$$\sum_{n=1}^\infty p_n(x,y) < \infty$$
for all $x,y \in I$. This statement can be proven in the same was as \cite[Theorem 2.6]{NT}. For the convenience of the reader, we give the argument. Denote by $\dimmin(y)$ the dimension of the carrier Hilbert space of $y$, when $y$ is viewed as an irreducible representation of $A_u(I_2)$. Since $F$ is not a multiple of a unitary $2 \times 2$ matrix, it follows that $\dimq(y) > \dimmin(y)$ for all $y \in I \setminus \{\eps\}$. Denote by $\mult(z;y_1 \ot \cdots \ot y_n)$ the multiplicity of the irreducible representation $z \in I$ in the tensor product of the irreducible representations $y_1,\ldots,y_n$. Since the fusion rules of $A_u(F)$ and $A_u(I_2)$ are identical, it follows that
$$\mult(z;y_1 \ot \cdots \ot y_n) \leq \dimmin(y_1) \cdots \dimmin(y_n) \; .$$
One then computes for all $x,y \in I$, $n \in \N$,
\begin{align*}
p_n(x,y) &= \sum_{z \subset \overline{x} \ot y} \mu^{*n}(z) \frac{\dimq(y)}{\dimq(x)\dimq(z)} \\
&= \frac{\dimq(y)}{\dimq(x)} \sum_{z \subset \overline{x} \ot y} \; \sum_{y_1,\ldots,y_n \in I} \mult(z;y_1 \ot \cdots y_n) \frac{\mu(y_1) \cdots \mu(y_n)}{\dimq(y_1) \cdots \dimq(y_n)} \\
&\leq \frac{\dimq(y)}{\dimq(x)} \dim(\overline{x} \ot y) \; \rho^n
\end{align*}
where $\displaystyle \rho = \sum_{y \in I} \mu(y) \frac{\dimmin(y)}{\dimq(y)}$. Since $\mu$ is non-trivial and $F$ is not a multiple of a $2 \times 2$ unitary matrix, we have $0 < \rho < 1$. Transience of the random walk follows immediately.

An element $x \in I$ is said to be \emph{indecomposable} if $x = y \ot z$ implies $y = \eps$ or $z = \eps$. Equivalently, $x$ is an alternating product of the letters $\alpha$ and $\beta$.

For every $x \in I$, we denote by $\dimq(x)$ the \emph{quantum dimension} of the irreducible representation labeled by $x$. Since $\dimq(\al) > 2$, take $0 < q < 1$ such that $\dimq(\alpha) = \dimq(\beta) = q + 1/q$. An important part of the proof of Theorem \ref{thm.main} is based on the technical estimates provided by Lemma \ref{lem.approx} and they require $q < 1$, i.e.\ $\dimq(\al) > 2$.

Denote the $q$-factorials $$[n]_q := \frac{q^n - q^{-n}}{q - q^{-1}} = q^{n-1} + q^{n-3} + \cdots + q^{-n+3} + q^{-n+1} \; .$$
Writing $x = x_1 \ot \cdots \ot x_n$ where the words $x_1,\ldots,x_n$ are indecomposable, we have
\begin{equation}\label{eq.quant-dim}
\dimq(x) = [|x_1|+1]_q \cdots [|x_n|+1]_q \; .
\end{equation}
For later use, note that it follows that
\begin{equation}\label{eq.quant-dim-ineq}
\dimq(xy) \geq q^{-|y|} \dimq(x)
\end{equation}
for all $x,y \in I$.

Whenever $x \in I \cup \partial I$, we denote by $[x]_n$ the word consisting of the first $n$ letters of $x$ and by $[x]^n$ the word that arises by removing the first $n$ letters from $x$. So, by definition, $x = [x]_n [x]^n$.

\section{Poisson boundary of the classical random walk on $\Irred(\cG)$}

Given a probability measure $\mu$ on $I:= \Irred(\cG)$, the Markov operator $P_\mu : \ell^\infty(\cGh) \recht \ell^\infty(\cGh)$ preserves the center $\cZ(\ell^\infty(\cGh)) = \ell^\infty(I)$ and hence, defines an ordinary random walk on the countable set $I$ with $n$-step transition probabilities
\begin{equation}\label{eq.transition}
p_x p_n(x,y) = p_x P_\mu^n(p_y) \; .
\end{equation}
As shown above, this random walk is transient whenever $\mu(\eps) < 1$.
Denote by $\Hcentr(\cGh,\mu)$ the commutative von Neumann algebra of bounded $P_\mu$-harmonic functions in $\ell^\infty(I)$, with product given by $a \cdot b = \lim_n P_\mu^n(ab)$ and the sequence being strongly$^*$-convergent. Write $p(x,y) = p_1(x,y)$.

The set $I$ becomes in a natural way a tree, the Cayley tree of the semi-group $\N * \N$. Whenever $\mu$ is a generating probability measure with finite support, the following properties of the random walk on $I$ can be checked easily.
\begin{itemize}
\item Uniform irreducibility: there exists an integer $M$ such that, for any pair $x,y\in I$ of neighboring edges of the tree, there exists an integer $k\leq M$, such that $p_k(x,y)>0$.
 \item Bounded step-length: there exists an integer $N$ such that $p(x,y)>0$ implies that $d(x,y)\leq N$ where $d(x,y)$ equals the length of the unique geodesic path from $x$ to $y$.
 \item There exists a $\delta>0$ such that $p(x,y)>0$ implies that $p(x,y)\geq \delta$.
\end{itemize}

So, we can apply \cite[Theorem 2]{PiWo} and identify the Poisson boundary of the random walk on $I$, with the boundary $\partial I$ of infinite words in $\al$,$\beta$, equipped with a probability measure in the following way.

\begin{theorem}[{Picardello and Woess, \cite[Theorem 2]{PiWo}}] \label{thm.classical}
Let $\mu$ be a finitely supported generating measure on $I = \Irred(A_u(F))$, where $F$ is not a scalar multiple of a $2 \times 2$ unitary matrix. Consider the associated random walk on $I$ with transition probabilities given by \eqref{eq.transition} and the compactification $I \cup \partial I$ of $I$.
\begin{itemize}
\item The random walk converges almost surely to a point in $\partial I$.
\item Denote, for every $x \in I$, by $\nu_x$ the hitting probability measure on $\partial I$, where $\nu_x(\cU)$ is defined as the probability that the random walk starting in $x$ converges to a point in $\cU$. Then, the formula
\begin{equation}\label{eq.classicalPI}
\Upsilon(F)(x) = \int_{\partial I} F (z) d\nu_x(z)
\end{equation}
defines a $*$-isomorphism $\Upsilon : \rL^\infty(\partial I,\nu_\eps) \recht \Hcentr(\cGh,\mu)$.
\end{itemize}
\end{theorem}

In fact, Theorem 2 in \cite{PiWo}, identifies $\partial I$ with the \emph{Martin compactification} of the given random walk on $I$. It is a general fact (see \cite[Theorem 24.10]{woess}), that a transient random walk converges almost surely to a point of the minimal Martin boundary and that the hitting probability measures provide a realization of the Poisson boundary through the Poisson integral formula \eqref{eq.classicalPI}, see \cite[Theorem 24.12]{woess}.

The rest of this section is devoted to the proof of the non-atomicity of the harmonic measures $\nu_x$.

\begin{lemma} \label{lem.cont}
For all $x,y\in I$ and $z\in \partial I$, the sequence
$$\Bigl(\frac{\dimq(x[z]_n)}{\dimq(y[z]_n)}\Bigr)_n$$
converges. By a slight abuse of notation, we denote the limit by $\dimq\bigl(\frac{xz}{yz} \bigr)$. The map $\partial I \recht \R_+ : z\mapsto \dimq\bigl(\frac{xz}{yz} \bigr)$ is continuous.
\end{lemma}
\begin{proof}
Fix $x,y \in I$. Whenever $z \in \partial I$ and $n \in \N$, denote
$$f_n(z) = \frac{\dimq(x[z]_n)}{\dimq(y[z]_n)} \; .$$
If $z \not\in \{\al \be \al \cdots, \be \al \be \cdots \}$, write $z = z_1 \ot z_2$ for some $z_1 \in I$, $z_1 \neq \eps$ and some $z_2 \in \partial I$. Denote by $\cU$ the neighborhood of $z$ consisting of words of the form $z_1z' = z_1 \ot z'$. For all $s \in \cU$ and all $n \geq |z_1|$, we have
$$f_n(s) = \frac{\dimq(x z_1)}{\dimq(y z_1)} \; .$$
Hence, for all $s \in \cU$, the sequence $n \mapsto f_n(s)$ is eventually constant and converges to a limit that is constant on $\cU$.

Also for $z \in \{\al \be \al \cdots, \be \al \be \cdots \}$, the sequence $f_n(z)$ is convergent. Take $z = \al \be \al \cdots$. Write $x = x_0 \ot x_1$ where $x_1$ is the longest possible (and maybe empty) indecomposable word ending with $\be$. Write $y = y_0 \ot y_1$ similarly. It follows that
$$f_n(z) = \frac{\dimq(x_0)}{\dimq(y_0)} \; \frac{[n+|x_1|+1]_q}{[n+|y_1|+1]_q} \recht \frac{\dimq(x_0)}{\dimq(y_0)} \; q^{|y_1|-|x_1|} \; .$$
The convergence of $f_n(z)$ for $z = \be \al \be \cdots$ is proven analogously.

Write $f(z) = \lim_n f_n(z)$. It remains to prove that $f$ is continuous in $z = \al \be \al \cdots$ and in $z = \be \al \be \cdots$. In both cases, define for every $n \in \N$, the neighborhood $\cU_n$ of $z$ consisting of all $s \in \partial I$ with $[s]_n = [z]_n$. For every $s \in \cU$, $s \neq z$, there exists $m \geq n$ such that $f(s) = f_m(z)$. The continuity of $f$ in $z$ follows.
\end{proof}

Whenever $x,y \in I \cup \partial I$, define $(x|y) := \max \{n \mid [x]_n = [y]_n\}$.

\begin{lemma} \label{lemma.formula-atoms}
For all $x \in I$ and all $w \in \partial I$, we have
$$\nu_x(w) = \frac{1}{\dimq(x)} \sum_{k=0}^{(x|w)} \dimq \left(\frac{[w]_k[w]^k}{\overline{[x]^k} [w]^k}\right)\nu_{\eps}\bigl(\overline{[x]^k} [w]^k \bigr) \; .$$
\end{lemma}
\begin{proof}
By Theorem \ref{thm.classical}, our random walk converges almost surely to a point of $\partial I$ and we denoted by $\nu_x$ the hitting probability measure. So, $(\psi_x \ot \psi_{\mu^{* n}})\deh \recht \nu_x$ weakly$^*$ in $\rC(I \cup \partial I)^*$.

Recall that $\cE : \ell^\infty(\cGh) \recht \ell^\infty(I)$ denotes the conditional expectation defined by $\cE(b) p_y = \psi_y(b) p_y$. Whenever $|z| > |x|$, we have
$$\cE \bigl( (\psi_x \ot \id)\deh(p_z)\bigr) = \sum_{k=0}^{(x|z)} \frac{\dimq(z)}{\dimq(x) \, \dimq\bigl(\overline{[x]^k} [z]^k \bigr)} p_{\overline{[x]^k} [z]^k} \; .$$
Denote $q_z = \sum_{s \in I} p_{zs}$ and observe that $q_z \in \rC(I \cup \partial I)$. It follows that for all $|z| > |x|$,
$$\cE \bigl( (\psi_x \ot \id)\deh(q_z)\bigr) = \sum_{k=0}^{(x|z)} f_k$$
where $f_k \in \ell^\infty(I)$ is defined by $f_k(y) = 0$ if $y$ does not start with $\overline{[x]^k} [z]^k$ and
$$f_k(\overline{[x]^k} [z]^k y) = \frac{1}{\dimq(x)} \frac{\dimq(z y)}{\dimq\bigl(\overline{[x]^k} [z]^k y\bigr)} \; .$$
By Lemma \ref{lem.cont}, we have $f_k \in \rC(I \cup \partial I)$.

Defining $\cU_z$ as the subset of $\partial I$ consisting of infinite words starting with $z$, it follows that, for $|z| > |x|$
\begin{equation}\label{eq.convolutie}
\nu_x(\cU_z) = \sum_{k=0}^{(x|z)} \int_{\partial I} f_k(y) \; d\nu_\eps (y) \; .
\end{equation}
The lemma follows by letting $z \recht w$.
\end{proof}

\begin{proposition} \label{prop.non-atomic}
The support of the harmonic measure $\nu_\eps$ is the whole of $\partial I$.
The harmonic measure $\nu_{\eps}$ has no atoms in words ending with $\al \be \al \be \cdots$
\end{proposition}

\begin{remark} \label{rem.non-atomic}
The same methods as in the proof of Proposition \ref{prop.non-atomic} given below, but involving more tedious computations, show in fact that $\nu_\eps$ is non-atomic. To prove our main theorem, it is only crucial that $\nu_\eps$ has no atoms in words ending with $\al \be \al \be \cdots$. We believe that it should be possible to give a more conceptual proof of the non-atomicity of $\nu_\eps$ and refer to \cite[Proposition 8.3.10]{VV-thesis} for an ad hoc proof along the lines of the proof of Proposition \ref{prop.non-atomic}.
\end{remark}

\begin{proof}[Proof of Proposition \ref{prop.non-atomic}]
Because of Lemma \ref{lemma.formula-atoms} and the equality
$$\nu_\eps  = \sum_{x \in I} \mu^{*k}(x) \nu_x$$
for all $k \geq 1$, we observe that if $w$ is an atom for $\nu_\eps $, then all $w'$ with the same tail as $w$ are atoms for all $\nu_x$, $x \in I$. Denote by $\cU_z \subset \partial I$ the open subset of $\partial I$ consisting of infinite words that start with $z$. A similar argument, using \eqref{eq.convolutie}, shows that if $\nu_\eps(\cU_z) = 0$, then $\nu_x(\cU_y) = 0$ for all $x,y \in I$, which is absurd because $\nu_x$ is a probability measure. So, the support of $\nu_\eps$ is the whole of $\partial I$.

So, we assume now that $w:= \al \be \al \be \cdots$ is an atom for $\nu_\eps$ and derive a contradiction.

Denote by $\delta_{w}$ the function on $\partial I$ that is equal to $1$ in $w$ and $0$ elsewhere. Using the $*$-isomorphism in Theorem \ref{thm.classical}, it follows that the bounded function
$$\xi \in \ell^\infty(\cGh) : \xi(x) := \nu_x(w) = \int_{\partial I} \delta_{w} \; d\nu_x$$
is harmonic.

We will prove that $\xi$ attains its maximum and apply the maximum principle for irreducible random walks (see e.g.\ \cite[Theorem 1.15]{woess}) to deduce that $\xi$ must be constant. This will lead to a contradiction.

Denote
$$w^\al_n := \underbrace{\al \be \al \cdots}_{n \; \text{letters}} \quad\text{and}\quad w^\be_n := \underbrace{\be \al \be \cdots}_{n \; \text{letters}} \; .$$
Note that all elements of $I$ are either of the form
$$w^\al_{2n+1} x \quad\text{where}\;\; n \in \N \;\;\text{and}\;\; x \in \{\eps\} \cup \al I$$
or of the form
$$w^\al_{2n} x \quad\text{where}\;\; n \in \N \;\;\text{and}\;\; x \in \{\eps\} \cup \be I \; .$$
By Lemma \ref{lemma.formula-atoms} and formula \eqref{eq.quant-dim}, we get that for $n \in \N$ and $x \in \{\eps\} \cup \al I$,
\begin{align*}
\xi(w^\al_{2n+1} x) & = \sum_{k=0}^{2n+1} \frac{1}{[2(n+1)]_q \; \dimq(x)^2} \; \dimq\Bigl(\frac{w^\al_k \hspace{1cm} [w]^k}{w^\be_{2n+1-k} \; [w]^k}\Bigr) \; \nu_\eps(\overline{x} \be \al \be \cdots) \\
& = \sum_{k=0}^{2n+1} \frac{1}{[2(n+1)]_q \; \dimq(x)^2} \; q^{2(n-k) +1} \; \nu_\eps(\overline{x} \be \al \be \cdots) = \frac{\nu_\eps(\overline{x} \be \al \be \cdots)}{\dimq(x)^2} \; .
\end{align*}
Since $\nu_\eps$ is a probability measure, it follows that $x \mapsto \xi(w^\al_{2n+1} x)$ is independent of $n$ and summable over the set $\{\eps\} \cup \al I$. Analogously, it follows that $x \mapsto \xi(w^\al_{2n} x)$ is independent of $n$ and summable over the set $\{\eps\} \cup \be I$. As a result, $\xi$ attains its maximum on $I$. By the maximum principle, $\xi$ is constant. Since $\xi(\eps) \neq 0$, this constant is non-zero and we arrive at a contradiction with the summability of $x \mapsto \xi(w^\al_{2n+1} x)$ over the infinite set $\{\eps\} \cup \al I$.
\end{proof}

\section{Topological boundary and boundary action for the dual of $A_u(F)$}\label{sec.boundary}

Before proving Theorem \ref{thm.main}, we construct a compactification for $\cGh$, i.e.\ a unital C$^*$-algebra $\mB$ lying between $\rc_0(\cGh)$ and $\ell^\infty(\cGh)$. This C$^*$-algebra $\mB$ is a non-commutative version of $\rC(I \cup \partial I)$. The construction of $\mB$ follows word by word the analogous construction given in \cite[Section 3]{VaVe} for $\cG = A_o(F)$. So, we only indicate the necessary modifications.

For all $x,y \in I$ and $z \subset x \ot y$, we choose an isometry $V(x \ot y,z) \in \Mor(x \ot y,z)$. Since $z$ appears with multiplicity one in $x \ot y$, the isometry $V(x \ot y,z)$ is uniquely determined up to multiplication by a scalar $\lambda \in S^1$. Therefore, the following unital completely positive maps are uniquely defined (cf.\ \cite[Definition 3.1]{VaVe}).

\begin{definition} \label{def.psi}
Let $x,y \in I$. We define unital completely positive
maps
$$\psi_{xy,x} : \cL(H_x) \recht \cL(H_{xy}) : \psi_{xy,x}(A) = V(x \ot y,xy)^* (A \ot 1) V(x \ot y,xy) \; .$$
\end{definition}

\begin{theorem} \label{thm.VV}
The maps $\psi_{xy,x}$ form an inductive system of completely positive maps. Defining
$$\mB=\{a\in \ell^\infty(\cGh)\ \mid\ \quad
 \forall \veps >0,\ \exists n\in \mathbb{N}\quad \mbox{such that} \quad\Vert ap_{xy}-\psi_{xy,x}(ap_x)\Vert<\veps$$ $$ \quad\mbox{for all}\quad x,y\in I\quad \mbox{with}\quad \lx\geq n\} \; ,$$
we get that $\mB$ is a unital C$^*$-subalgebra of $\ell^\infty(\cGh)$ containing $\rc_0(\cGh)$.
\begin{itemize}
\item The restriction of the comultiplication $\deh$ yields a left action $\be_{\cGh}$ of $\cGh$ on $\mB$~:
$$\be_{\cGh} : \mB \recht \M(\rc_0(\cGh) \ot \mB) : a \mapsto \deh(a) \; .$$
\item The restriction of the adjoint action  of $\cG$ on $\ell^\infty(\cGh)$ yields a right action of $\cG$ on $\mB$~:
$$\be_{\cG} : \mB \recht \mB \ot \rC(\cG) : a \mapsto \bV(a \ot 1) \bV^*\; .$$
Here, $\bV \in \ell^\infty(\cGh) \ovt \rL^\infty(\cG)$ is defined as $\bV = \sum_{x \in I} U^x$. The action $\be_\cG$ is continuous in the sense that $\lspan \be_\cG(\mB) (1 \ot \rC(\cG))$ is dense in $\mB \ot \rC(\cG)$.
\end{itemize}
\end{theorem}
\begin{proof}
One can repeat word by word the proofs of \cite[Propositions 3.4 and 3.6]{VaVe}. The crucial ingredients of these proofs are the approximate commutation formulae provided by \cite[Lemmas A.1 and A.2]{VaVe} and they have to be replaced by the inequalities provided by Lemma \ref{lem.approx}.
\end{proof}

We denote $\mB_\infty := \mB/\rc_0(\cGh)$ and call it the topological boundary of $\cGh$. Both actions $\be_{\cG}$ and $\be_{\cGh}$ preserve the ideal $\rc_0(\cGh)$ and hence yield actions on $\mB_\infty$ that we still denote by $\be_{\cG}$ and $\be_{\cGh}$.

We identify $\rC(I \cup \partial I)$ with $\mB \cap \cZ(\ell^\infty(\cGh)) = \mB \cap \ell^\infty(I)$. Similarly, $\rC(\partial I) \subset \mB_\infty$.

We partially order $I$ by writing $x \leq y$ if $y = x z$ for some $z \in I$. Define
$$\psi_{\infty,x} : \cL(H_x) \recht \mB : \psi_{\infty,x}(A)p_y = \begin{cases} \psi_{y,x}(A) \;\;\text{if}\;\; y \geq x \\ 0 \;\;\text{else}\end{cases} \; .$$
We use the same notation for the composition of $\psi_{\infty,x}$ with the quotient map $\mB \recht \mB_\infty$, yielding the map $\psi_{\infty,x} : \cL(H_x) \recht \mB_\infty$.

\begin{lemma}\label{lem.field}
The inclusion $\rC(\partial I) \subset \mB_\infty$ defines a continuous field of unital C$^*$-algebras. Denote, for every $x \in \partial I$, by $J_x$ the closed two-sided ideal of $\mB_\infty$ generated by the functions in $\rC(\partial I)$ vanishing in $x$.

For every $x = x_1 \ot x_2 \ot \cdots$ in $\partial_0 I$, there exists a unique surjective $*$-homomorphism $$\pi_x : \mB_\infty \recht \bigotimes_{k=1}^\infty \cL(H_{x_k})$$
satisfying $\Ker \pi_x = J_x$ and $\pi_x(\psi_{\infty,x_1\cdots x_n}(A)) = A \ot 1$ for all $A \in \bigotimes_{k=1}^n \cL(H_{x_k}) = \cL(H_{x_1\cdots x_n})$.
\end{lemma}
\begin{proof}
Given $x \in \partial I$, define the decreasing sequence of projections $e_n \in \mB$ given by
$$e_n := \sum_{y \in I} p_{[x]_n y} \; .$$
Denote by $\pi : \mB \recht \mB_\infty$ the quotient map. It follows that
\begin{equation}\label{eq.norm}
\|\pi(a) + J_x\| = \lim_n \|a e_n \|
\end{equation}
for all $a \in \mB$.

To prove that $\rC(\partial I) \subset \mB_\infty$ is a continuous field, let $y \in I$, $A \in \cL(H_y)$ and define $a \in \mB$ by $a:= \psi_{\infty,y}(A)$. Put $f : \partial I \recht \R_+ : f(x) = \|\pi(a) + J_x\|$. We have to prove that $f$ is a continuous function. Define $\cU \subset \partial I$ consisting of infinite words starting with $y$. Then, $\cU$ is open and closed and $f$ is zero, in particular continuous, on the complement of $\cU$. Assume that the last letter of $y$ is $\al$ (the other case, of course, being analogous). If $x \in \cU$ and $x \neq y \be \al \be \al \cdots$, write $x = yz \ot u$ for some $z \in I$, $u \in \partial I$. Define $\cV$ as the neighborhood of $x$ consisting of infinite words of the form $yzu'$ where $u' \in \partial I$ and $yzu' = yz \ot u'$. Then, $f$ is constantly equal to $\|\psi_{yz,y}(A)\|$ on $\cV$. It remains to prove that $f$ is continuous in $x := y \be \al \be \al \cdots$. Let
$$w_n = \underbrace{\be \al \be \cdots}_{n \; \text{letters}} \; .$$
Then, the sequence $\|\psi_{yw_n,y}(A)\|$ is decreasing and converges to $f(x)$. If $\cU_n$ is the neighborhood of $x$ consisting of words starting with $yw_n$, it follows that
$$f(x) \leq f(u) \leq \|\psi_{yw_n,y}(A) \|$$
for all $u \in \cU_n$. This proves the continuity of $f$ in $x$. So, $\rC(\partial I) \subset \mB_\infty$ is a continuous field of C$^*$-algebras.

Let now $x = x_1 \ot x_2 \ot \cdots$ be an element of $\partial_0 I$. Put $y_n = x_1 \ot \cdots \ot x_n$ and
$$f_n := \sum_{z \in I} p_{y_n z} \; .$$
The map $A \mapsto f_{n+1} \psi_{\infty,y_n}(A)$ defines a unital $*$-homomorphism from $\cL(H_{y_n})$ to $f_{n+1} \mB$. Since $1 - f_{n+1} \in J_x$, we obtain the unital $*$-homomorphism $\theta_n : \cL(H_{y_n}) \recht \mB/J_x$. The $*$-homomorphisms $\theta_n$ are compatible and combine into the unital $*$-homomorphism
$$\theta : \bigotimes_{k=1}^\infty \cL(H_{x_k}) \recht \mB/J_x \; .$$
By \eqref{eq.norm}, $\theta$ is isometric. Since the union of all $\psi_{\infty,y_n}(\cL(H_{y_n})) + J_x$, $n \in \N$, is dense in $\mB$, it follows that $\theta$ is surjective. The composition of the quotient map $\mB \recht \mB/J_x$ and the inverse of $\theta$ provides the required $*$-homomorphism $\pi_x$.
\end{proof}

\section{Proof of Theorem \ref{thm.main}}

We prove Theorem \ref{thm.main} by performing the following steps.
\begin{itemize}
\item Construct on the boundary $\mB_\infty$ of $\cGh$, a faithful KMS state $\om_\infty$, to be considered as the harmonic state and satisfying $(\psi_\mu \ot \om_\infty)\be_{\cGh} = \om_\infty$. Extend $\be_{\cGh}$ to an action
    $$\be_{\cGh} : (\mB_\infty,\om_\infty)\dpr \recht \ell^\infty(\cGh) \ovt (\mB_\infty,\om_\infty)\dpr$$
    and denote by $\PI := (\id \ot \om_\infty) \be_{\cGh}$ the \emph{Poisson integral.}
\item Prove a quantum \emph{Dirichlet property~:} for all $a \in \mB$, we have $\PI(a) - a \in \rc_0(\cGh)$. It will follow that $\PI$ is a normal and faithful $*$-homomorphism of $(\mB_\infty,\om_\infty)\dpr$ onto a von Neumann subalgebra of $\poisson$.
\item By Theorem \ref{thm.classical}, $\PI$ is a $*$-isomorphism of $\rL^\infty(\partial I,\nu_\eps ) \subset (\mB_\infty,\om_\infty)\dpr$ onto $\Hcentr(\cGh,\mu)$. Deduce that the image of $\PI$ is the whole of $\poisson$.
\item Use Lemma \ref{lem.field} to identify $(\mB_\infty,\om_\infty)\dpr$ with a field of ITPFI factors.
\end{itemize}

\begin{proposition}
The sequence $\psi_{\mu^{*n}}$ of states on $\mB$ converges weakly$^*$ to a KMS state $\om_\infty$ on $\mB$. The state $\om_\infty$ vanishes on $\rc_0(\cGh)$. We still denote by $\om_\infty$ the resulting KMS state on $\mB_\infty$. Then, $\om_\infty$ is faithful on $\mB_\infty$.

We have $(\psi_\mu \ot \om_\infty)\be_{\cGh} = \om_\infty$, so that we can uniquely extend $\be_{\cGh}$ to an action
$$\be_{\cGh} : (\mB_\infty,\om_\infty)\dpr \recht \ell^\infty(\cGh) \ovt (\mB_\infty,\om_\infty)\dpr$$
that we still denote by $\be_{\cGh}$.

The state $\om_\infty$ is invariant under the action $\be_\cG$ of $\cG$ on $\mB_\infty$. We extend $\be_\cG$ to an action on $(\mB_\infty,\om_\infty)\dpr$ that we still denote by $\be_\cG$.

The normal, completely positive map
\begin{equation}\label{eq.PI}
\PI : (\mB_\infty,\om_\infty)\dpr \recht \poisson : \PI = (\id \ot \om_\infty)\be_{\cGh}
\end{equation}
is called the \emph{Poisson integral.} It satisfies the following properties (recall that $\al_{\cGh}$ and $\al_\cG$ were introduced in Definition \ref{def.poisson}).
\begin{itemize}
\item $\epsh \circ \PI = \om_\infty$.
\item $(\PI \ot \id)\circ \be_\cG = \al_\cG \circ \PI$.
\item $(\id \ot \PI) \circ \be_{\cGh} = \al_{\cGh} \circ \PI$.
\end{itemize}
For every $x = x_1 \ot x_2 \ot \cdots$ in $\partial_0 I$, denote by $\om_x$ the infinite tensor product state on $\bigotimes_{k=1}^\infty \cL(H_{x_k})$, of the states $\psi_{x_k}$ on $\cL(H_{x_k})$. Using the notation $\pi_x$ of Lemma \ref{lem.field}, we have
\begin{equation}\label{eq.productstate}
\om_\infty(a) = \int_{\partial_0 I} \om_x(\pi_x(a)) \; d\nu_\eps(x)
\end{equation}
for all $a \in \mB_\infty$.
\end{proposition}
\begin{proof}
Define the $1$-parameter group of automorphisms $(\si_t)_{t \in \R}$ of $\ell^\infty(\cGh)$ given by
$$\si_t(a) p_x = Q_x^{it} a p_x Q_x^{-it} \; .$$
Since $\si_t(\psi_{\infty,x}(A)) = \psi_{\infty,x}(Q_x^{it} A Q_x^{-it})$, it follows that $(\si_t)$ is norm-continuous on the C$^*$-algebra $\mB$.

By Theorem \ref{thm.classical}, the sequence of probability measures $\mu^{* n}$ on $I \cup \partial I$ converges weakly$^*$ to $\nu_\eps$. It follows that $\psi_{\mu^{* n}}(a) \recht 0$ whenever $a \in \rc_0(\cGh)$.
Given $x \in I$ and $A \in \cL(H_x)$, put $a := \psi_{\infty,x}(A)$. Denote by $\partial(xI)$ the set of infinite words starting with $x$ and by $\partial_0(xI)$ its intersection with $\partial_0(I)$.
We get, using Proposition \ref{prop.non-atomic},
\begin{align*}
\psi_{\mu^{* n}}(a) &= \sum_{y \in I} \mu^{*n}(y) \psi_y(\psi_{\infty,x}(A)) = \sum_{y \in xI} \mu^{*n}(y) \psi_x(A) \\ & \recht \psi_x(A) \nu_\eps(\partial(xI)) = \psi_x(A) \nu_\eps(\partial_0(xI)) = \int_{\partial_0 I} \om_y(\pi_y(a)) \; d\nu_\eps(y) \; .
\end{align*}
So, the sequence $\psi_{\mu^{*n}}$ of states on $\mB$ converges weakly$^*$ to a state on $\mB$ that we denote by $\om_\infty$ and that satisfies \eqref{eq.productstate}. Since all
$\psi_{\mu^{*n}}$ satisfy the KMS condition w.r.t.\ $(\si_t)$, also $\om_\infty$ is a KMS state. If $a \in \mB_\infty^+$ and $\om_\infty(a) = 0$, it follows from \eqref{eq.productstate} that $\om_x(\pi_x(a)) = 0$ for $\nu_\eps$-almost every $x \in \partial_0 I$. Since $\om_x$ is faithful, it follows that $\| \pi(a) + J_x \| = 0$ for $\nu_\eps$-almost every $x \in \partial I$. By Proposition \ref{prop.non-atomic}, the support of $\nu_\eps$ is the whole of $\partial I$ and by Lemma \ref{lem.field}, $x \mapsto \|\pi(a) + J_x\|$ is a continuous function. It follows that $\|\pi(a) + J_x\| = 0$ for all $x \in \partial I$ and hence, $a = 0$. So, $\om_\infty$ is faithful.

Since $(\psi_\mu \ot \psi_{\mu^{*n}})\be_{\cGh} = \psi_{\mu^{*(n+1)}}$, it follows that $(\psi_\mu \ot \om_\infty)\be_{\cGh} = \om_\infty$. So, $(\psi_{\mu^{*k}} \ot \om_\infty)\be_{\cGh} = \om_\infty$ for all $k \in \N$. Since $\mu$ is generating, there exists for every $x \in I$, a $C_x > 0$ such that $(\psi_x \ot \om_\infty)\be_{\cGh} \leq C_x \om_\infty$. As a result, we can uniquely extend $\be_{\cGh}$ to a normal $*$-homomorphism
$$(\mB_\infty,\om_\infty)\dpr \recht \ell^\infty(\cGh) \ovt (\mB_\infty,\om_\infty)\dpr \; .$$
Since $\be_{\cGh}$ is an action, the same holds for the extension to the von Neumann algebra $(\mB_\infty,\om_\infty)\dpr$.

Because $(\psi_\mu \ot \om_\infty)\be_{\cGh} = \be_{\cGh}$ and because $\be_{\cGh}$ is an action, the Poisson integral defined by \eqref{eq.PI} takes values in $\poisson$. It is straightforward to check that $\PI$ intertwines $\be_\cG$ with $\al_\cG$ and $\be_{\cGh}$ with $\al_{\cGh}$.
\end{proof}

\begin{theorem}\label{thm.quantumdirichlet}
The compactification $\mB$ of $\cGh$ satisfies the quantum Dirichlet property, meaning that, for all $a \in \mB$,
\begin{equation*}
\| (\PI(a) - a)p_x \| \to 0
\end{equation*}
if $\lx\to\infty$.

In particular, the Poisson integral $\PI$ is a normal and faithful $*$-homomorphism of $(\mB_\infty,\om_\infty)\dpr$ onto a von Neumann subalgebra of $\poisson$.
\end{theorem}

We will deduce Theorem \ref{thm.quantumdirichlet} from the following lemma.

\begin{lemma} \label{lem.dirichlet}
For every $a\in \mB$ we have that
\begin{equation}\label{eq.dirichlet}
\sup_{y \in I} \| (\id \ot \psi_y)\deh(a)p_x - ap_x \| \recht 0
\end{equation}
when $\lx \recht \infty$.
\end{lemma}

\begin{proof}
Fix $a \in \mB$ with $\|a\| \leq 1$.
Choose $\veps > 0$. Take $n$ such that $\|a p_{x_0x_1} - \psi_{x_0x_1,x_0}(a p_{x_0})\| < \veps$ for all $x_0,x_1 \in I$ with $|x_0| = n$.

Denote $d_{S^1}(V,W) = \inf \{ \|V - \lambda W\| \mid \lambda \in S^1\}$. By \eqref{eq.variant} below, take $k$ such that
\begin{equation}\label{eq.useful}
\begin{split}
d_{S^1} \bigl( (V(x_0 \ot x_1x_2, x_0x_1x_2) \ot 1) & V(x_0x_1x_2 \ot \overline{x_2} u, x_0 x_1 u) \; , \\ & (1 \ot V(x_1 x_2 \ot \overline{x_2} u,x_1 u))V(x_0 \ot x_1 u,x_0x_1 u) \bigr) < \frac{\veps}{2}
\end{split}
\end{equation}
whenever $|x_1| \geq k$.

Finally, take $l$ such that $q^{2l} < \veps$. We prove that
\begin{equation}\label{eq.5eps}
\| (\id \ot \psi_y)\deh(a)p_x - ap_x \| <  5 \veps
\end{equation}
for all $x,y \in I$ with $|x| \geq n+k+l$.

Choose $x,y \in I$ with $|x| \geq n+k+l$ and write $x = x_0 x_1 x_2$ with $|x_0| = n$, $|x_1|=k$ and hence, $|x_2| \geq l$. We get
\begin{align*}
(\id \ot \psi_y)\deh(a) p_x &= \sum_{z \subset x \ot y} (\id \ot \psi_y)\bigl( V(x \ot y,z) a p_z V(x \ot y,z)^* \bigr) \\ & = \sum_{z \subset x \ot y} \frac{\dimq(z)}{\dimq(x) \; \dimq(y)} V(z \ot \overline{y},x)^* (a p_z \ot 1) V(z \ot \overline{y},x) \\
& = \sum_{z \subset x_2 \ot y} \frac{\dimq(x_0x_1z)}{\dimq(x) \; \dimq(y)} V(x_0x_1z \ot \overline{y},x)^* (a p_{x_0x_1z} \ot 1) V(x_0x_1 z \ot \overline{y},x) \\ & \hspace{1cm} + \sum \text{remaining terms} \; .
\end{align*}
In order to have remaining terms, $y$ should be of the form $y = \overline{x_2} y_0$ and then, using \eqref{eq.quant-dim-ineq} and the assumption $\|a\| \leq 1$,
\begin{align*}
\sum \|\text{remaining terms}\| &= \sum_{z \subset x_0x_1 \ot y_0} \frac{\dimq(z)}{\dimq(x_0x_1x_2) \; \dimq(\overline{x_2} y_0)} \\ & \leq
\sum_{z \subset x_0x_1 \ot y_0} q^{2|x_2|} \frac{\dimq(z)}{\dimq(x_0x_1) \; \dimq(y_0)} = q^{2|x_2|} < \veps \; .
\end{align*}
Combining this estimate with the fact that $\|a p_{x_0x_1z} - \psi_{x_0x_1z,x_0}(a p_{x_0})\| < \veps$, it follows that
\begin{align*}
\bigl\| & (\id \ot \psi_y)\deh(a) p_x - a p_x \bigr\| \\ & \leq 2 \veps + \Bigl\| a p_x - \sum_{z \subset x_2 \ot y} \frac{\dimq(x_0x_1z)}{\dimq(x) \; \dimq(y)} V(x_0x_1z \ot \overline{y},x)^* ( \psi_{x_0x_1z,x_0}(a p_{x_0}) \ot 1) V(x_0x_1 z \ot \overline{y},x) \Bigr\| \; .
\end{align*}
But now, \eqref{eq.useful} implies that
$$\bigl\| (\id \ot \psi_y)\deh(a) p_x - a p_x \bigr\| \leq 3 \veps + \Bigl\| a p_x - \sum_{z \subset x_2 \ot y} \frac{\dimq(x_0x_1z)}{\dimq(x) \; \dimq(y)} \psi_{x,x_0}(a p_{x_0}) \Bigr\| \; .$$
Since $\|\psi_{x,x_0}(a p_{x_0}) - a p_x \| < \veps$ and $\|a \| \leq 1$, we get
$$\bigl\| (\id \ot \psi_y)\deh(a) p_x \bigr\| \leq 4 \veps + \Bigl(1 - \sum_{z \subset x_2 \ot y} \frac{\dimq(x_0x_1z)}{\dimq(x) \; \dimq(y)} \Bigr) \; .$$
The second term on the right hand side is zero, unless $y = \overline{x_2} y_0$, in which case, it equals
$$\sum_{z \subset x_0x_1 \ot y_0} \frac{\dimq(z)}{\dimq(x_0x_1x_2) \; \dimq(\overline{x_2} y_0)} \leq \sum_{z \subset x_0x_1 \ot y_0} q^{2|x_2|} \frac{\dimq(z)}{\dimq(x_0x_1) \; \dimq(y_0)} \leq \veps$$
because of \eqref{eq.quant-dim-ineq}. Finally, \eqref{eq.5eps} follows and the lemma is proven.
\end{proof}

\begin{proof}[Proof of Theorem \ref{thm.quantumdirichlet}]
Let $a \in \mB$. Given $\veps > 0$, Lemma \ref{lem.dirichlet} provides $k$ such that
$$\| (\id \ot \psi_{\mu^{* n}})\deh(a)p_x - ap_x \| \leq \veps$$
for all $n \in \N$ and all $x$ with $|x| \geq k$. Since $\psi_{\mu^{* n}}\to \om_\infty$ weakly$^*$, it follows that
$$
\| (\PI(a) - a)p_x \| \leq \eps
$$
whenever $|x| \geq k$. This proves \eqref{eq.dirichlet}.

It remains to prove the multiplicativity of $\PI$. We know that $\PI:\mB_\infty\to \poisson$ is a unital, completely positive map. Since $\epsh \circ \PI = \om_\infty$, $\PI$ is faithful. Denote by $\pi:\poisson\to \frac{\ell^\infty(\cGh)}{\rc_0(\cGh)}$ the quotient map, which is also a unital, completely positive map. By \eqref{eq.dirichlet}, we have $\pi \circ \PI = \id$. So, for all $a \in \mB_\infty$, we find
$$\pi(\PI(a)^* \cdot \PI(a)) \leq \pi(\PI(a^* a)) = a^* a = \pi(\PI(a))^* \pi(\PI(a)) \leq \pi(\PI(a)^* \cdot \PI(a))$$
We claim that $\pi$ is faithful. If $a \in \poisson^+ \cap \rc_0(\cGh)$, we have $\epsh(a) = \psi_{\mu^{*n}}(a)$ for all $n$ and the transience of $\mu$ combined with the assumption $a \in \rc_0(\cGh)$, implies that $\epsh(a) = 0$ and hence, $a = 0$. So, we conclude that $\PI(a)^* \cdot \PI(a) = \PI(a^* a)$ for all $a \in \mB_\infty$. Hence, $\PI$ is multiplicative on $\mB_\infty$ and also on $(\mB_\infty,\om_\infty)\dpr$ by normality.
\end{proof}

\begin{remark}\label{rem.reinterpretation}
We now give a reinterpretation of Theorem \ref{thm.classical}. Denote by $\Omega = I^\N$ the path space of the random walk with transition probabilities \eqref{eq.transition}. Elements of $\Omega$ are denoted by $\ux$, $\underline{y}$, etc. For every $x \in I$, one defines the probability measure $\bP_x$ on $\Omega$ such that $\bP_x(\{x\} \times I \times I \times \cdots) = 1$ and
$$\bP_x(\{(x,x_1,x_2,\cdots,x_n)\} \times I \times I \times \cdots) = p(x,x_1) \; p(x_1,x_2) \; \cdots \; p(x_{n-1},x_n) \; .$$
Choose a probability measure $\eta$ on $I$ with $I = \supp \eta$. Write $\bP = \sum_{x \in I} \eta(x) \bP_x$.

Define on $\Omega$ the following equivalence relation: $\ux \sim \underline{y}$ iff there exist $k,l \in \N$ such that $x_{n+k} = y_{n+l}$ for all $n \in \N$. Whenever $F \in \Hcentr(\cGh,\mu)$, the martingale convergence theorem implies that the sequence of measurable functions $\Omega \recht \C : \ux \mapsto F(x_n)$ converges $\bP$-almost everywhere to a $\sim$-invariant bounded measurable function on $\Omega$, that we denote by $\pi_\infty(F)$.
Denote by $\rL^\infty(\Omega/_\sim,\bP)$ the von Neumann subalgebra of $\sim$-invariant functions in $\rL^\infty(\Omega,\bP)$.
As such, $\pi_\infty : \Hcentr(\cGh,\mu) \recht \rL^\infty(\Omega/_\sim,\bP)$ is a $*$-isomorphism.

By Theorem \ref{thm.quantumdirichlet}, we can define the measurable function $\bnd : \Omega \recht \partial I$ such that $\bnd \ux = \lim_n x_n$ for $\bP$-almost every $\ux \in \Omega$ and where the convergence is understood in the compact space $I \cup \partial I$. Recall that we denote, for $x \in I$, by $\nu_x$ the hitting probability measure on $\partial I$. So, $\nu_x(A) = \bP_x(\bnd^{-1}(A))$ for all measurable $A \subset \partial I$ and all $x \in I$.

Again by Theorem \ref{thm.quantumdirichlet}, $\pi_\infty \circ \Upsilon$ is a $*$-isomorphism of $\rL^\infty(\partial I, \nu_\eps)$ onto 
$\rL^\infty(\Omega/_\sim,\bP)$. We claim that for all $F \in \rL^\infty(\partial I, \nu_\eps)$, we have
$$\bigl((\pi_\infty \circ \Upsilon)(F)\bigr)(\ux) = F(\bnd \ux) \quad\text{for $\bP$-almost every}\;\; \ux \in \Omega \;  .$$
Let $A \subset \partial I$ be measurable. Define $F_n : \Omega \recht \R : F_n(\ux) = \nu_{x_n}(A)$. Then, $F_n$ converges almost everywhere with limit equal to $(\pi_\infty \circ \Upsilon)(\chi_A)$. If the measurable function $G : \Omega \recht \C$ only depends on $x_0,\ldots,x_k$, one checks that
$$\int_{\Omega} F_n(\ux) \; G(\ux) \; d\bP(\ux) = \int_{\bnd^{-1}(A)} G(\ux) \; d\bP(\ux) \quad\text{for all}\;\; n > k \; .$$
From this, the claim follows.

Since the $*$-isomorphism $\pi_\infty \circ \Upsilon$ is given by $\bnd$, it follows that for every $\sim$-invariant bounded measurable function $F$ on $\Omega$, there exists a bounded measurable function $F_1$ on $\partial I$ such that $F(\ux) = F_1(\bnd \ux)$ for $\bP$-almost every path $\ux \in \Omega$.
\end{remark}

We are now ready to prove the main Theorem \ref{thm.main}.

\begin{proof}[Proof of Theorem \ref{thm.main}]
Because of Theorem \ref{thm.quantumdirichlet} and Lemma \ref{lem.field}, it remains to show that
$$\PI : (\mB_\infty,\om_\infty)\dpr \recht \poisson$$
is surjective.

Whenever $\gamma : N \recht N \ovt \rL^\infty(\cG)$ is an action of $\cG$ on the von Neumann algebra $N$, we denote, for $x \in I$, by $N^x \subset N$ the \emph{spectral subspace} of the irreducible representation $x$. By definition, $N^x$ is the linear span of all $S(H_x)$, where $S$ ranges over the linear maps $S : H_x \recht N$ satisfying $\gamma(S(\xi)) = (S \ot \id)(U_x (\xi \ot 1))$. The linear span of all $N^x$, $x \in I$, is a weakly dense $*$-subalgebra of $N$, called the spectral subalgebra of $N$. For $n \in \N$, we denote by $N^n$ the linear span of all $N^x$, $|x| \leq n$.

Fixing $x,y \in I$, consider the adjoint action $\gamma : \cL(H_{xy}) \recht \cL(H_{xy}) \ot \rC(\cG)$ given by $\gamma(A) = U_{xy}(A \ot 1) U_{xy}^*$. The fusion rules of $\cG = A_u(F)$ imply that $\cL(H_{xy})^{2|x|} = \psi_{xy,x}(\cL(H_x))$.

For the rest of the proof, put $M := (\mB_\infty,\om_\infty)\dpr$. We use the action $\be_\cG$ of $\cG$ on $M$ and the action $\al_\cG$ of $\cG$ on $\poisson$. Fix $k \in \N$. It suffices to prove that $\poisson^k \subset \PI(M)$.

Define, for all $y\in I$, the subset
\begin{align*}
 V_y:=\{yz\mid z\in I \;\;\text{and}\;\; yz = y \ot z \} \; .
\end{align*}
Define the projections $$q_y=\sum_{z\in V_y}p_{z}\in\mB$$
and consider $q_y$ also as an element of the von Neumann algebra $M$. Define $W_y \subset \partial I$ as the subset of infinite words of the form $yu$, where $u \in \partial I$ and $yu = y \ot u$.

Define the unital $*$-homomorphism
$$
\zeta:\rC(W_y) \ot \cL(H_y)\to M q_y: \begin{cases} \zeta(F)p_{yz}=\psi_{yz,y}(F(yz)) & \quad \text{when}\;\; y\ot z=yz\\  0 & \quad \text{else}\end{cases} \; .
$$
{\bf Claim.} For all $y\in I$, there exists a linear map
\begin{equation*}
T_y :\poisson^{2|y|}\cdot \PI(q_y) \subset \poisson \to \rL^\infty(W_y)\ot \cL(H_y)
\end{equation*}
satisfying the following conditions.
\begin{itemize}
\item $T_y$ is isometric for the $2$-norm on $\poisson$ given by the state $\epsh$ and the $2$-norm on $\rL^\infty(W_y)\ot \cL(H_y)^k$ given by the state $\nu_\veps\ot \psi_y$.
\item $(T_y \circ \PI \circ \zeta)(F) = F$ for all $F \in \rC(W_y)\ot \cL(H_y)$.
\end{itemize}
To prove this claim, we use the notations and results introduced in Remark \ref{rem.reinterpretation}. Fix $y \in I$. Consider $a\in \poisson^{2|y|}\cdot \PI(q_y)$. If $\underline{x}\in \Omega$ is such that $\operatorname{bnd}(\underline{x})\in W_y$, then, for $n$ big enough, $x_n$ will be of the form $x_n=y\ot z_n$. By definition of $\alpha_\cG$, we have that $a p_{x_n}\in \cL(H_{x_n})^{2|y|}$. So, we can take elements $a_{\underline{x},n}\in \cL(H_y)$ such that $a p_{x_n}=\psi_{x_n,y}(a_{\underline{x},n})$. We prove that, for $\bP$-almost every path $\underline{x}$ with  $\bnd \underline{x} \in W_y$, the sequence $(a_{\underline{x},n})_n$ is convergent. We then define $T_y(a) \in \rL^\infty(W_y) \ot \cL(H_y)$ such that $T_y(a)(\bnd \ux) = \lim_n a_{\ux,n}$ for $\bP$-almost every path $\ux$ with $\bnd \ux \in W_y$.

Take $d\in \cL(H_y)$. Then, for $\bP$-almost every path $\underline{x}$ such that $\bnd \underline{x} \in W_y$ and $n$ big enough, we get that
$$
\psi_y(da_{\underline{x},n}) =\psi_{x_n}\bigl(\psi_{x_n,y}(da_{\underline{x},n}) \bigr)=\psi_{x_n}\bigl(\psi_{x_n,y}(d)\psi_{x_n,y}(a_{\underline{x},n}) \bigr) =\psi_{x_n}\bigl(\psi_{x_n,y}(d)a p_{x_n} \bigr)
$$
In the second step, we used the multiplicativity of $\psi_{x_n,y}:\cL(H_y)\to \cL(H_{x_n})$ which follows because $x_n=y\ot z_n$. Also note that $\Vert a_{\underline{x},n}\Vert\leq \Vert a\Vert$.
From Theorem \ref{thm.quantumdirichlet}, it follows that
\begin{equation*}
\bigl\Vert\PI\bigl(\zeta(1 \ot d)\bigr)p_{x_n}-\psi_{x_n,y}(d)p_{x_n}\bigr\Vert\to 0
\end{equation*}
whenever $x_n$ converges to a point in $W_y$. This implies that
$$\bigl\vert\psi_y(da_{\underline{x},n})- \psi_{x_n}\bigl(\PI\bigl(\zeta(1 \ot d)\bigr) a p_{x_n}\bigr)\bigr\vert\to 0$$
for $\bP$-almost every path $\underline{x}$ with $\bnd \underline{x} \in W_y$.

From \cite[Proposition 3.3]{INT}, we know that for $\bP$-almost every path $\ux$,
$$\big|\psi_{x_n}\bigl( \PI(\zeta(1 \ot d)) a p_{x_n} \bigr)p_{x_n} - \cE\bigl( \PI(\zeta(1 \ot d)) \cdot a \bigr) p_{x_n} \big| \recht 0 \; .$$
As before, $\cE(b)p_x = \psi_x(b) p_x$. It follows that
$$\bigl\vert\psi_y(da_{\underline{x},n}) p_{x_n} - \cE\bigl(\PI\bigl(\zeta(1 \ot d)\bigr)\cdot a \bigr) p_{x_n}\vert\to 0 \; .$$

Note that $\cE$ maps $\poisson$ onto $\Hcentr(\cGh,\mu)$. Whenever $F \in \Hcentr(\cGh,\mu)$, the sequence $F(x_n)$ converges for $\bP$-almost every path $\ux$. We conclude that for every $d \in \cL(H_y)$, the sequence $\psi_y(da_{\underline{x},n})$ is convergent for $\bP$-almost every path $\ux$ with $\bnd \underline{x} \in W_y$. Since $\cL(H_y)$ is finite dimensional, it follows that the sequence $(a_{\ux,n})_n$ in $\cL(H_y)$ is convergent for $\bP$-almost every path $\underline{x}$ with $\bnd \underline{x}\in W_y$.

By Remark \ref{rem.reinterpretation}, we get $T_y(a) \in \rL^\infty(W_y) \ot \cL(H_y)$ such that $T_y(a)(\bnd \ux) = \lim_n a_{\ux,n}$ for $\bP$-almost every path $\underline{x}$ with $\bnd \underline{x}\in W_y$. From the definition of $a_{\ux,n}$, we get that
\begin{equation}\label{eq.def}
\bigl\Vert\psi_{x_n,y}\bigl(T_y(a)(\bnd \underline{x} )\bigr)-a p_{x_n}\bigr\Vert\to 0
\end{equation}
for $\bP$-almost every path $\underline{x}$ such that $\bnd \underline{x} \in W_y$.

The map $T_y$ is isometric. Indeed, by the defining property \eqref{eq.def} and again by \cite[Proposition 3.3]{INT}, we have, for $\bP$-almost every path $\ux$ with $\bnd \ux \in W_y$,
$$\psi_y\bigl( T_y(a)(\bnd \underline{x} )^*T_y(a)(\bnd \underline{x} ) \bigr) = \lim_{n \recht \infty} \psi_{x_n}(a^* a p_{x_n}) = (\pi_\infty \circ \cE)(a^* \cdot a)(\ux) \; .$$
Here, $\pi_\infty$ denotes the $*$-isomorphism $\Hcentr(\cGh,\mu) \recht \rL^\infty(\Omega/_\sim,\bP)$ introduced in Remark \ref{rem.reinterpretation}. On the other hand, by Remark \ref{rem.reinterpretation}, $\bigl((\pi_\infty \circ \PI)(q_y)\bigr)(\ux) = 0$ for $\bP$-almost every path $\ux$ with $\bnd \ux \not\in W_y$. Since $$\int_{\Omega} \bigl((\pi_\infty \circ \cE)(b)\bigr)(\ux) \; d\bP_\eps(\ux) = \epsh(b)$$
for all $b \in \poisson$, it follows that $T_y$ is an isometry in $2$-norm.

We next prove that $(T_y \circ \PI \circ \zeta)(F) = F$ for all $F \in \rC(W_y)\ot \cL(H_y)$. Let $\atil \in \rC(I \cup \partial I) \subset \ell^\infty(I)$ and let $a$ be the restriction of $\atil$ to $\partial I$. Take $A \in \cL(H_y)$. It suffices to take $F = a \ot A$. Theorem \ref{thm.quantumdirichlet} implies that
$$
\bigl\Vert \atil p_{x_n} \psi_{x_n,y}(A) -(\PI \circ \zeta)(a \ot A)p_{x_n}\bigr\Vert\to 0
$$
for $\bP$-almost every path $\underline{x}$. On the other hand, for $\bP$-almost every path $\ux$ with $\bnd \ux \in W_y$, the scalar $\atil p_{x_n}$ converges to $a(\bnd \ux)$. In combination with \eqref{eq.def}, it follows that $(T_y \circ \PI \circ \zeta)(a \ot A) = a \ot A$, concluding the proof of the claim.

Having proven the claim, we now show that for all $y \in I$, $\poisson^{2|y|} \cdot \PI(q_y) \subset \PI(M)$. Take $a \in \poisson^{2|y|} \cdot \PI(q_y)$. Let $d_n$ be a bounded sequence in the C$^*$-algebra $\rC(W_y) \ot \cL(H_y)$ converging to $T_y(a)$ in $2$-norm. Since $T_y \circ \PI$ is an isometry in $2$-norm, it follows that $\zeta(d_n)$ is a bounded sequence in $M$ that converges in $2$-norm. Denoting by $c \in M$ the limit of $\zeta(d_n)$, we conclude that $T_y(\PI(c)) = T_y(a)$ and hence, $\PI(c) = a$.

Fix $k \in \N$. A fortiori, $\poisson^k\cdot \PI(q_y) \subset \PI(M)$ for all $y \in I$ with $2|y| \geq k$. By Proposition \ref{prop.non-atomic}, the harmonic measure $\nu_\eps$ has no atoms in infinite words ending with $\al \be \al \be \cdots$. As a result, $1$ is the smallest projection in $M$ that dominates all $q_y$, $y \in I$, $2|y| \geq k$. So, $\poisson^k \subset \PI(M)$ for all $k \in \N$. This finally implies that $\PI$ is surjective.
\end{proof}

\section{Solidity and the Akemann-Ostrand property}

In Section \ref{sec.boundary}, we followed the approach of \cite{VaVe} to construct the compactification $\mB$ of $\cGh$. In fact, more of the constructions and results of \cite{VaVe} carry over immediately to the case $\cG = A_u(F)$. We continue to assume that $F$ is not a multiple of a $2 \times 2$ unitary matrix.

Denote by $\rL^2(\cG)$ the GNS Hilbert space defined by the Haar state $h$ on $\rC(\cG)$. Denote by $\lambda : \rC(\cG) \recht \cL(\rL^2(\cG))$ the corresponding GNS representation and define $\Cred(\cG) := \lambda(\rC(\cG))$. We can view $\lambda$ as the left-regular representation. We also have a right-regular representation $\rho$ and the operators $\lambda(a)$ and $\rho(b)$ commute for all $a,b \in \rC(\cG)$ (see \cite[Formulae (1.3)]{VaVe}).

Repeating the proofs of \cite[Proposition 3.8 and Theorem 4.5]{VaVe}, we arrive at the following result.

\begin{theorem}\label{thm.AO}
The boundary action $\be_{\cGh}$ of $\cGh$ on $\mB$ defined in Theorem \ref{thm.VV} is
\begin{itemize}
\item \emph{amenable} in the sense of \cite[Definition 4.1]{VaVe};
\item \emph{small at infinity} in the following sense: the comultiplication $\deh$ restricts as well to a right action of $\cGh$ on $\mB$; this action leaves $\rc_0(\cGh)$ globally invariant and becomes the trivial action on the quotient $\mB_\infty$.
\end{itemize}
\end{theorem}

By construction, $\mB$ is a nuclear C$^*$-algebra and hence, as in \cite[Corollary 4.7]{VaVe}, we get that
\begin{itemize}
\item $\cG$ satisfies the Akemann-Ostrand property: the homomorphism
$$\Cred(\cG) \otalg \Cred(\cG) \recht \frac{\cL(\rL^2(\cG))}{\cK(\rL^2(\cG))} : a \ot b \mapsto \lambda(a)\rho(b) + \cK(\rL^2(\cG))$$
is continuous for the minimal C$^*$-tensor product $\otmin$.
\item $\Cred(\cG)$ is an exact C$^*$-algebra.
\end{itemize}

As before, we denote by $\rL^\infty(\cG)$ the von Neumann algebra acting on $\rL^2(\cG)$ generated by $\lambda(\rC(\cG))$. From
\cite[Th\'{e}or\`{e}me 3]{Ba2}, it follows that $\rL^\infty(\cG)$ is a factor, of type II$_1$ if $F$ is a multiple of an $n \times n$ unitary matrix and of type III in the other cases.

Applying \cite[Theorem 6]{ozawa} (in fact, its slight generalization provided by \cite[Theorem 2.5]{VaVe}), we get the following corollary of Theorem \ref{thm.AO}. Recall that a II$_1$ factor $M$ is called \emph{solid} if for every diffuse von Neumann subalgebra $A \subset M$, the relative commutant $M \cap A'$ is injective. An arbitrary von Neumann algebra $M$ is called \emph{generalized solid} if the same holds for every diffuse von Neumann subalgebra $A \subset M$ which is the image of a faithful normal conditional expectation.

\begin{corollary}
When $n \geq 3$ and $\cG = A_u(I_n)$, the II$_1$ factor $\rL^\infty(\cG)$ is solid. When $n \geq 2$, $F \in \GL(n,\C)$ is not a multiple of an $n \times n$ unitary matrix and $\cG = A_u(F)$, the type III factor $\rL^\infty(\cG)$ is generalized solid.
\end{corollary}

\section{Appendix: approximate intertwining relations}

We fix an invertible matrix $F$ and assume that $F$ is not a scalar multiple of a unitary $2 \times 2$ matrix. Define $\cG = A_u(F)$ and label the irreducible representations of $\cG$ by the monoid $\N * \N$, freely generated by $\al$ and $\be$. The representation labeled by $\al$ is the fundamental representation of $\cG$ and $\be$ is its contragredient. Define $0 < q < 1$ such that $\dimq(\al) = \dimq(\be) = q + \frac{1}{q}$. Recall from Section \ref{sec.boundary} that whenever $z \subset x \ot y$, we choose an isometry $V(x \ot y,z) \in \Mor(x \ot y,z)$. Observe that $V(x \ot y,z)$ is uniquely determined up to multiplication by a scalar $\lambda \in S^1$. We denote by $p^{x \ot y}_z$ the projection $V(x \ot y,z)V(x \ot y,z)^*$.

\begin{lemma}\label{lem.approx}
There exists a constant $C > 0$ that only depends on $q$ such that
\begin{align}
\| \bigl( V(xr \ot \overline{r} y, xy) \ot 1_z \bigr)p^{xy \ot z}_{xyz} - (1_{xr} \ot p^{\overline{r}y \ot z}_{\overline{r}yz})\bigl(V(xr \ot \overline{r}y,xy) \ot 1_z\bigr)\| & \leq C q^{|y|} \; ,\notag\\
\| \bigl( 1_x \ot V(yr \ot \overline{r}z,yz)\bigr) p^{x \ot yz}_{xyz} - (p^{x \ot yr}_{xyr} \ot 1_{\overline{r}z}) \bigl( 1_x \ot V(yr \ot \overline{r}z,yz)\bigr)\| & \leq C q^{|y|}\label{eq.estimate2}
\end{align}
for all $x,y,z,r \in I$.
\end{lemma}

One way of proving Lemma \ref{lem.approx} consists in repeating step by step the proof of \cite[Lemma A.1]{VaVe}. But, as we explain now, Lemma \ref{lem.approx} can also be deduced more directly from \cite[Lemma A.1]{VaVe}.

\begin{proof}[Sketch of proof]
Whenever $y = y_1 \ot y_2$ with $y_1 \neq \eps \neq y_2$, the expressions above are easily seen to be $0$. Denote
$$v_n = \underbrace{\al \ot \be \ot \al \ot \cdots}_{n \;\;\text{tensor factors}} \quad\text{and}\quad w_n = \underbrace{\be \ot \al \ot \be \ot \cdots}_{n \;\;\text{tensor factors}} \; .$$
The remaining estimates that have to be proven, reduce to estimates of norms of operators in $\Mor(v_n,v_m)$ and $\Mor(w_n,w_m)$. Putting these spaces together in an infinite matrix, one defines the C$^*$-algebras
$$A := \bigl(\Mor(v_n,v_m)\bigr)_{n,m} \quad\text{and}\quad B := \bigl(\Mor(w_n,w_m)\bigr)_{n,m}$$
generated by the subspaces $\Mor(v_n,v_m)$ and $\Mor(w_n,w_m)$, respectively.
Choose unit vectors $t \in \Mor(\al \ot \be,\eps)$ and $s \in \Mor(\be \ot \al,\eps)$ such that $(t^* \ot 1)(1 \ot s) = (q+1/q)^{-1}$. By \cite[Lemme 5]{Ba2}, the C$^*$-algebra $A$ is generated by the elements $1^{\ot 2k} \ot t \ot 1^{\ot l}$, $1^{\ot (2k+1)} \ot s \ot 1^{\ot l}$. A similar statement holds for $B$.

Denote by $U$ the fundamental representation of the quantum group $\SU_{-q}(2)$ and let $t_0 \in \Mor(U \ot U,\eps)$ be a unit vector. The proofs of \cite[Theorems 5.3 and 6.2]{BDRV} (which heavily rely on results in \cite{Ba1,Ba2}), imply the existence of $*$-isomorphisms
$$\pi_A : (\Mor(U^{\ot n} , U^{\ot m}))_{n,m} \recht A \quad\text{and}\quad \pi_B : (\Mor(U^{\ot n} , U^{\ot m}))_{n,m} \recht B$$
satisfying
$$\pi_A(1^{\ot 2k} \ot t_0 \ot 1^{\ot l})=1^{\ot 2k} \ot t \ot 1^{\ot l} \quad\text{and}\quad \pi_A(1^{\ot (2k+1)} \ot t_0 \ot 1^{\ot l}) = 1^{\ot (2k+1)} \ot s \ot 1^{\ot l}$$
and similarly for $\pi_B$.

As a result, the estimates to be proven, follow directly from the corresponding estimates for $\SU_{-q}(2)$ proven in \cite[Lemma A.1]{VaVe}.
\end{proof}

Using the notation
$$d_{S^1}(V,W) = \inf \{ \| V - \lambda W \| \mid \lambda \in S^1 \} \; ,$$
several approximate commutation relations can be deduced from Lemma \ref{lem.approx}. For instance, after a possible increase of the constant $C$, \eqref{eq.estimate2} implies that
\begin{equation}\label{eq.variant}
d_{S^1}\Bigl( \bigl(1_x \ot V(yr \ot \overline{r}z,yz)\bigr) V(x \ot yz,xyz) \; , \; \bigl( V(x \ot yr,xyr) \ot 1_{\overline{r}z}\bigr) V(xyr \ot \overline{r}z,xyz) \Bigr) \leq C q^{|y|}
\end{equation}
for all $x,y,z,r \in I$. We again refer to \cite[Lemma A.1]{VaVe} for a full list of approximate intertwining relations.

\end{document}